\documentclass{ifacconf}

\usepackage{graphicx}      
\usepackage{natbib}        
\usepackage{amsmath}
\usepackage{amssymb}
\usepackage{mathdots}
\usepackage{natbib}
\usepackage[ruled,vlined,algo2e]{algorithm2e}
\usepackage[all,cmtip]{xy}
\usepackage{tikz}
\usetikzlibrary{automata,positioning,shapes,arrows,shadows}       

\usepackage{xcolor} 

\begin{document}
\begin{frontmatter}

\title{Computing Truncated Joint Approximate Eigenbases for Model Order Reduction} 

\thanks[footnoteinfo]{Loring acknowledges partial support from the National Science Foundation \#2110398. Vides acknowledges partial support from the Scientific Computing Innovation Center of UNAH under project PI-174-DICIHT.}

\author[Terry]{Terry Loring} 
\author[Fredy]{Fredy Vides} 

\address[Terry]{Department of Mathematics and Statistics, University of New Mexico, 
   Albuquerque, (e-mail: loring@math.unm.edu).}
\address[Fredy]{Scientific Computing Innovation Center, School of Mathematics and Computer Science, 
   Universidad Nacional Aut\'onoma de Honduras, Tegucigalpa (e-mail: fredy.vides@unah.edu.hn)}

\end{frontmatter}

\section{Introduction}
Consider a collection of $d$ Hermitian matrices $X_1,\ldots,X_d$ in $\mathbb{R}^{n\times n}$ and a $d$-tuple $\boldsymbol{\lambda}=(\lambda_1,\ldots,\lambda_d)\in \mathbb{R}^d$. Let us consider the problem determined by the computation of a collection of joint approximate eigenvectors that can be represented as a rectangular matrix $W \in \mathbb{C}^{n\times r}$ with orthonormal columns such that
\begin{equation}
W=\arg\min_{\hat{W}\in \mathbb{R}^{n \times r}}\sum_{j=1}^d \left\|X_j\hat{W}-\hat{W}\Lambda_j\right\|_F^2.
\label{eq:main_partial_iso}
\end{equation}
Solutions to problem \eqref{eq:main_partial_iso} can be used for model order reduction as will be illustrated in \S\ref{sec:experiments}.

Given one Hermitian matrix $X$ we are only interested in
the real part of the pseudospectrum. By the usual definition, real
$\lambda$ is in the $\epsilon$-pseudospectrum of $X$ if
\[
\left\Vert (X-\lambda)^{-1}\right\Vert ^{-1}\leq\epsilon.
\]
One can easily see this is equivalent to the condition
\[
\exists\boldsymbol{v}\text{ such that }\|\boldsymbol{v}\|=1\text{ and }\left\Vert X\boldsymbol{v}-\lambda\boldsymbol{v}\right\Vert \leq\epsilon.
\]
We will call $\left\Vert X\boldsymbol{v}-\lambda\boldsymbol{v}\right\Vert $
the \emph{eigen-error}. This comes up all the time in applications,
and the less matrices commute the more it must be considered.

For Hermitian matrices $X_{1},X_{2},\dots,X_{d}$ we often want a
unit vector with the various eigen-errors small. There are many ways
to combine $d$ errors, such as their sum or maximum. Not surprisingly,
a clean theory arises when we consider the quadratic mean of the eigen-errors.

Here then is a definition of a pseudospectrum. In the noncommutative
setting, there are several notions of joint spectrum and joint pseudospectrum
that compete for our attention, such as one using Clifford algebras \citep{LoringPseudospectra}. None is best is all settings. 
\begin{defn}
Suppose we have finitely many Hermitian matrices $X_{1},X_{2}$ $,$ $\dots$ $,$ $X_{d}$.
Suppose $\epsilon>0$. A $d$-tuple $\boldsymbol{\lambda}$ is an
element of the quadratic $\epsilon$-pseudospectrum of $(X_{1},X_{2}$ $,$ $\dots$ $,$ $X_{d})$
if there exists as unit vector $\boldsymbol{v}$ so that
\begin{equation}
\sqrt{\sum_{j=1}^{d}\left\Vert X_{j}\boldsymbol{v}-\lambda_{j}\boldsymbol{v}\right\Vert ^{2}}\leq\epsilon.\label{eq:QMerrors}
\end{equation}
If (\ref{eq:QMerrors}) is true for $\epsilon=0$ then we say $\boldsymbol{\lambda}$
is an element of the quadratic spectrum of $(X_{1},X_{2}$ $,$ $\dots$ $,$ $X_{d})$.
The notation for the quadratic $\epsilon$-pseudospectrum of $(X_{1},X_{2}$ $,$ $\dots$ $,$ $X_{d})$
is $Q\Lambda_{\epsilon}(X_{1},X_{2},\dots,X_{d})$.
\end{defn}

\begin{rem}
Very simple examples show that the quadratic spectrum can often be
empty.

It should be said that the more interesting examples of this tend
to require calculation, or at least approximation, by numerical methods.
Often the best way to display the data is via images of 2D slices
through the function
\end{rem}

\[
\boldsymbol{\lambda}\mapsto\mu_{\boldsymbol{\lambda}}^{Q}(X_{1},\dots,X_{d})
\]
where we define 
\begin{equation}
\mu_{\boldsymbol{\lambda}}^{Q}(X_{1},\dots,X_{d})=\min_{\left\Vert \boldsymbol{v}\right\Vert =1}\sqrt{\sum_{j=1}^{d}\left\Vert X_{j}\boldsymbol{v}-\lambda_{j}\boldsymbol{v}\right\Vert ^{2}}.\label{eq:QuadPSfunction}
\end{equation}
That is, we have a measure of how good of a joint approximate eigenvector
we can find at $\boldsymbol{\lambda}$. Then, of course, the more
traditional interpretation of $Q\Lambda_{\epsilon}(X_{1},X_{2},\dots,X_{d})$ as the
sublevel sets of this function. 
\begin{rem}

We will make frequent use of the following notation:
\[
Q_{\boldsymbol{\lambda}}(X_{1},\dots,X_{d})=\sum_{j=1}^{d}\left(X_{j}-\lambda_{j}\right)^{2},
\]
\[
M_{\boldsymbol{\lambda}}(X_{1},\dots,X_{d})=\left[\begin{array}{c}
X_{1}-\lambda_{1}\\
\vdots\\
X_{d}-\lambda_{d}
\end{array}\right]
\]
Finally we use $\sigma_{\min}$ to indicate the smallest singular
value of a matrix. 
\end{rem}

As a particular application of quadratic pseudospectrum based techniques, for the computation of truncated joint approximate eigenbases, in section \S\ref{sec:experiments} we will present an application of these quadratic pseudospectral based methods to the computation of a reduced order model for a discrete-time system related to least squares realization of linear time invariant models \citep{DeMoor}.

\section{Main Results}
\label{sec:main}

We now list the main results that corresponding to some important properties of the quadratic pseudospectrum.

\begin{prop}
Suppose that $X_{1},X_{2},\dots,X_{d}$ are Hermitian matrices, that
$\epsilon>0$ and $\boldsymbol{\lambda}$ is in $\mathbb{R}^{d}$.
The following are equivalent.
\begin{enumerate}
\item $\boldsymbol{\lambda}$ is an element of the quadratic $\epsilon$-pseudospectrum
of $(X_{1},X_{2},\dots,X_{d})$;
\item $\sigma_{\min}\left(M_{\boldsymbol{\lambda}}(X_{1},\dots,X_{d})\right)\leq\epsilon$;
\item $\sigma_{\min}\left(Q_{\boldsymbol{\lambda}}(X_{1},\dots,X_{d})\right)\leq\epsilon^{2}$.
\end{enumerate}
\end{prop}

The following technical result is very helpful for numerical calculations.
Assuming that one does not care about the exact value of $\mu_{\boldsymbol{\lambda}}^{Q}(X_{1},\dots,X_{d})$
once this value is above some cutoff, then knowing Lipschitz continuity
allows one to skip calculating this values at many points near where
a high value has been found.
\begin{prop}
Suppose that $X_{1},X_{2},\dots,X_{d}$ are Hermitian matrices. The
function 
\[
\boldsymbol{\lambda}\mapsto\mu_{\boldsymbol{\lambda}}^{Q}(X_{1},\dots,X_{d}),
\]
with domain $\mathbb{R}^{d}$, is Lipschitz with Lipschitz constant
$1$.
\end{prop}

For details on the proofs of Propositions 4  and 5, the reader is kindly referred to \citep{arxiv.2204.10450}.

\section{Algorithm}
\label{sec:alg}
Combining the ideas and methods presented in \citep{7053905} and \citep{SC-siam}, with the ideas and results presented in \S\ref{sec:main}, we obtained Algorithm \ref{alg:pjad}.

\begin{algorithm2e}
\caption{Approximate Joint Eigenvectors Computation}
\label{alg:pjad}
\SetAlgoLined
\KwData{{\sc Hermitian matrices: $X_1,\ldots,X_d\in \mathbb{R}^{n\times n}$}, {\sc $d$-tuple $\lambda\in \mathbb{C}^d$}, {\sc Integer:} $1\leq k\leq n$, {\sc Threshold:} $\delta>0$, {\sc Selector:} $\phi$}
\KwResult{{\sc Partial isometry} $V\in \mathbb{O}(n,k)$}
\begin{itemize}
\item[0:] Set the choice indicator value $\phi$: $\phi=0$ for smallest eigenvalues or $\phi=1$ for largest eigenvales\;
\item[1:] Set $L:=\sum_{j=1}^N (X_j-\lambda_j I_n)^2$\;
\item[2:] Approximately solve $LV=V\Lambda$ for $V\in \mathbb{C}^{n\times k}, \Lambda\in \mathbb{C}^{k\times k}$ according to the flag value $\phi$\;
\For{$j\gets 1$ \KwTo $d$}{
    \begin{itemize}
    \item[3.0:] Set $Y_j:=V^\top (X_j-\lambda_j I_n) V$\;
    \item[3.1:] Set $Y_j:=(Y_j+Y_j^\top)/2$\;
    \end{itemize}
    }
\item[4:] Solve $W=\arg\min_{U\in \mathbb{O}(n)}\sum_{k=1}^d \mathrm{off}(U^\top Y_kU)$ using complex valued Jacobi-like techniques as in \cite{SC-siam} with threshold$=\delta$.\;
\item[5:] Set $V:=VW$\;
\end{itemize}
\KwRet{$V$}
\end{algorithm2e}

In this document, the operation $A^\top$ represents the transpose of some given matrix $A$.

\section{Example}
\label{sec:experiments}

Consider the discrete-time system with states $x_1(t)$ and $x_2(t)$ in $\mathbb{R}^{400}$:
\begin{align}
x_1(t+1) &= A_1x_1(t), \: x_2(t+1) = A_2x_1(t+1),\label{eq:LTI}\\
y_1(t) & = \hat{e}_{1,400}^\top x_1(t), \: y_2(t)  = \hat{e}_{2,400}^\top x_2(t),\nonumber
\end{align}
for some given matrices $A_1,A_2\in \mathbb{R}^{400\times 400}$ such that $A_1A_2=A_2A_1$ that are generated with the program {\tt QLMORDemo.py} available at \citep{FVides_TJAE}., here $\hat{e}_{1,400}$ and $\hat{e}_{2,400}$ denote the first and second columns of the identity matrix in $\mathbb{R}^{400\times 400}$, respectively. Let us consider the matrices 
\begin{align*}
H_1 &= A_1^\top A_1, \\
H_2 &= A_2^\top A_2, \\
H_3 &= A_1^\top A_2+A_2^\top A_1
\end{align*}
We can apply Algorithm \ref{alg:pjad} to $H_1,H_2,H_3$ with $\delta=10^{-5}$ obtaining the matrix $V\in \mathbb{R}^{400\times 6}$ with orthonormal columns, that can be used to compute a model order reduction for \eqref{eq:LTI}, determined by the following equations.
\begin{align*}
\hat{x}_1(t+1) &= V^\top A_1V \hat{x}_1(t), \: \hat{x}_2(t+1) = V^\top A_2V \hat{x}_1(t+1),\\
\hat{y}_1(t) & = \hat{e}_{1,400}^\top V\hat{x}_1(t), \: \hat{y}_2(t)  = \hat{e}_{2,400}^\top V\hat{x}_2(t).
\end{align*}

\begin{figure}
\centering
\includegraphics[scale=.6]{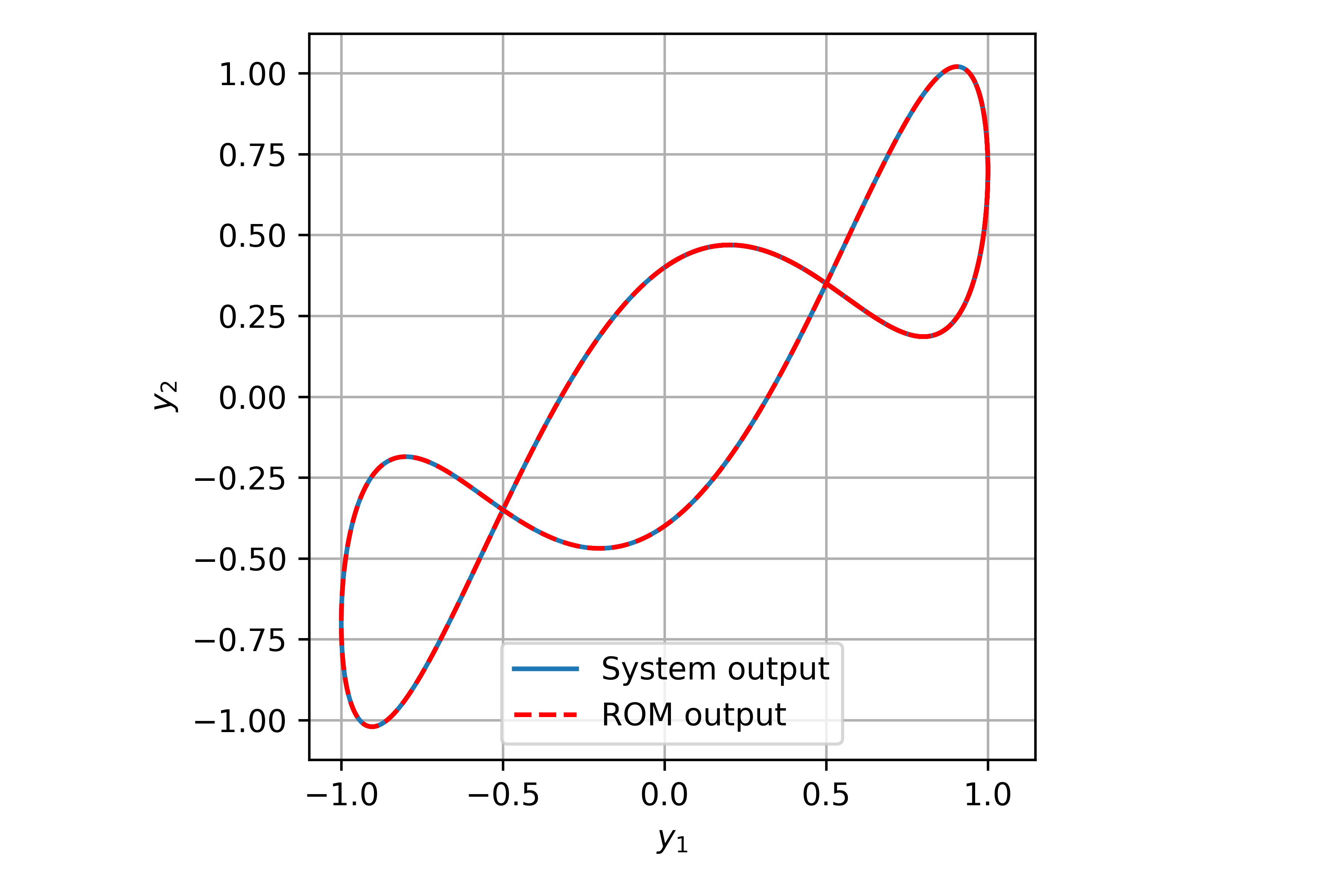}
\caption{Original system and ROM outputs.}
\label{fig:experiment_1_1}
\end{figure}

The outputs corresponding to the original and reduced order models are plotted in Figure \ref{fig:experiment_1_1}.

\bibliography{root}
\end{document}